\documentclass[11pt,a4paper,twoside]{amsart}

\usepackage{amsmath,amssymb,amsthm}
\usepackage{inputenc}
\usepackage{tikz} 
\usepackage{float}
\floatstyle{boxed} 
\usepackage{graphicx,subfigure}
\usepackage{caption}
\usepackage[T1]{fontenc} 
\usepackage{enumitem}
\usepackage{mathtools}
\usepackage{mathrsfs}
\usepackage{xcolor}
\usepackage[colorlinks=true, linkcolor=magenta, citecolor=cyan, urlcolor=blue,hyperfootnotes=false]{hyperref}
\usepackage{cleveref}

\usepackage[style=numeric,firstinits=true,sorting=nyt,backend=biber,maxnames=99]{biblatex}
\usepackage{etoolbox}
\usepackage{upgreek}

\def\@begintheorem#1#2{\par\bgroup{\sc #1 \ #2. }  \it \\\ignorespace }
\def\@opargbegintheorem#1#2#3{\par\bgroup{\sc #1\ #2 \ (#3).}  \it  \ignorespace}
\def\@endtheorem{\egroup}

\theoremstyle{plain}
\newtheorem{theorem}{Theorem}[section]

\theoremstyle{definition}

\theoremstyle{remark}
\newtheorem{remark}[theorem]{Remark}

\theoremstyle{plain}

\theoremstyle{plain}
\newtheorem{lemma}[theorem]{Lemma}

\theoremstyle{plain}

\theoremstyle{plain}
\newtheorem{proposition}[theorem]{Proposition}

\numberwithin{equation}{section}

\newcommand{\zb}{{\bar{z}}}
\newcommand{\C}{{\mathbb C}}

\newcommand{\N}{{\mathbb N}}
\newcommand{\D}{\mathcal{D}}

\newcommand{\R}{{\mathbb R}}

\newcommand{\V}{\mathbb{V}}

\newcommand{\Rt}{{\R}^3}
\newcommand{\Y}{{\mathcal{Y}}}
\newcommand{\hx}{\hat{x}}

\DeclareMathOperator{\Realpart}{Re}
\renewcommand{\Re}{\Realpart}

{\left\lbrace\begin{array}{@{}l@{}}}%
{\end{array}\right.}

\DeclarePairedDelimiter{\abs}{\lvert}{\rvert}
\DeclarePairedDelimiter{\norm}{\lVert}{\rVert}

\newcommand{\E}{\textbf{E}}

\newcommand{\ignora}[1]{}

\title{Sharp exponential localization for solutions of the Perturbed Dirac Equation}
\date{\today}
\author[B. Cassano]{Biagio Cassano}
\thanks{The author was partially supported by 
the ERC Advanced Grant 669689 HADE (European Research Council), by the MINECO project MTM2014-53145-P, by the Basque Government through the BERC 2014-2017 program, by the Spanish Ministry of Economy and Competitiveness MINECO: BCAM Severo Ochoa accreditation SEV-2013-0323,
by Istituto Italiano Di Alta Matematica, and
by the Czech Science Foundation (GA\v{C}R), grant No.~17-01706S.
}
\address{B.~Cassano, Department of Mathematics, %
  Universit\`{a} degli  Studi di Bari, %
  via Edoardo Orabona 4, 70125, Bari, Italy}
\email{biagio.cassano@uniba.it}

\usepackage{csquotes}

\bibliography{biblio.bib}
\begin{document}
\begin{abstract}
We determine the largest non-trivial rate of exponential decay at
infinity for solutions to the Dirac equation
 \begin{equation*}
   \D_n \psi + \V \psi = 0 
   \quad \text{ in }\R^n,  
 \end{equation*}
being
$\D_n$ the massless Dirac operator in dimension $n\geq 2$ and $\V$ a
(possibly non-Hermitian)
matrix-valued perturbation such that
$\abs{\V(x)} \sim \abs{x}^{-\epsilon}$ at infinity, for 
$-\infty < \epsilon < 1$.
Moreover, we show that our results are sharp
for $n =2,3$, providing explicit
examples of solutions that have the prescripted
decay, in presence of a potential with the related behaviour at
infinity.
\end{abstract}
\maketitle



\section{Introduction}
We investigate the rate of exponential decay at
infinity for solutions of the Dirac equation
 \begin{equation}\label{eq:main}
   \D_n \psi + \V \psi = 0 
   \quad \text{ in }\R^n,  
 \end{equation}
 where $\D_n$ is the massless Dirac operator in dimension $n \geq 2$,
 and $\V$ is a
(possibly non-Hermitian) matrix-valued perturbation such that
$\abs{\V(x)} \sim \abs{x}^{-\epsilon}$ at infinity, with 
$-\infty <\epsilon < -1$. 
We remark that if $\epsilon \leq 0$
we can also describe eigenfunctions of the massive Dirac operator
associated to any complex eigenvalue.

We prove a rigidity property of the Dirac operator: if $\psi$ is a solution to \eqref{eq:main} and 
has a too large exponential decay at infinity, then $\psi$ has
compact support.
Moreover, when $n=2,3$ we provide explicit examples of solutions to \eqref{eq:main}
that have the prescripted exponential decay at infinity, when the potential has
the related decay at infinity: this ensures that in these cases our results are sharp.
When informations on the local behaviour of $\V$ are
provided, our results determine also the sharpest
exponential decay at infinity  for non-trivial solutions
to \eqref{eq:main}.


Since the Dirac operator is the
matricial square root of the Laplace operator,
it is interesting to review the results on the analogous problem for
the latter. 
Let $u$ be a solution to
\begin{equation}\label{eq:introduzione1}
\Delta u + V u = Eu 
\quad \text{ in }\R^n,
\end{equation}
with $n \geq 2$, $E \in \R$, and assume that for some $\epsilon \in \R$
we have $\abs{V(x)} \leq C \abs{x}^{-\epsilon}$ for some $C>0$ and for big $\abs{x}$.
In \cite{froese19832} it is shown that if $\epsilon \geq 1/2$ and
$\exp[{\tau\abs{x}}] u \in L^2(\R^n) $ for a big enough $\tau \gg 1$,
then $u$ has compact support. 
In \cite{kondratiev1988qualitative}, E.~M.~Landis conjectured that
such phenomenon is more general, claiming that if
$V$ is bounded ($\epsilon = 0$) at infinity
and $\exp[{\tau\abs{x}}] u \in L^2(\R^n)$ for a big enough $\tau\gg 1$, 
then $u$ must have compact support.
This was disproved by Meshkov in \cite{meshkov1992possible}:
by means of the appropriate Carleman estimate he proved that 
a general solution $u$ to \eqref{eq:introduzione1} in $\R^n$, $n\geq 2$, has compact support if 
$\exp[{\tau\abs{x}^{4/3}}] u \in L^2(\R^n)$ for a big enough $\tau \gg
1$ and $V$ is bounded, i.e.~$\epsilon=0$.
Moreover he provided    an intelligent example in $\R^2$ of a bounded potential $V_M$ and a non-trivial
function $u_M(x) \sim \exp[-C \abs{x}^{4/3}]$ at infinity such that
\eqref{eq:introduzione1} holds true.
It is important to underline that the functions $u_M$ and $V_M$ are
\emph{complex valued}: in fact  the Carleman estimates can not
distinguish between real or complex valued functions,
and the question of the validity of Landis' conjecture was left open
in the case of real-valued potentials $V$.
The results in \cite{meshkov1992possible} were generalized to the
case $-\infty < \epsilon < 1/2$ in \cite{cruz1999unique}, where Cruz-Sampedro
showed that $u$ has compact support if 
$\exp[\tau\abs{x}^{p}] u \in L^2(\R^n) $ for $\tau \gg 1$ big enough
and $p=(4-2\epsilon)/3$,
and  provided  an example in $\R^2$ of a non-trivial function with the optimal decay, in presence
of a potential with the associated behaviour at infinity, in the style
of \cite{meshkov1992possible}.
In \cite{duyckaerts2008optimality} 
 Duyckaerts, Zuazua and Zhang observed 
that an easier construction for the examples with the critical decay
can be done considering \emph{vector-valued}
solutions to \eqref{eq:introduzione1}: 
a $\C^4$--valued solution is
constructed for \eqref{eq:introduzione1} in $\R^3$ and when $\epsilon
= 0$, but the matrix-valued potential $V$ has a logaritmic growth at infinity, 
that is almost optimal in this context.

In \cite{escauriaza2011unique} Escauriaza, Kenig, Ponce and Vega proved
unique continuation properties for solutions of the evolution 
Schr\"odinger equation with a time dependent potential
\begin{equation}\label{eq:introduzioneEKPV}
  i \partial_t u + \Delta u + V(x,t) u=0, 
  \quad
  (x,t) \in \R^n \times (0, +\infty),
\end{equation}
where $V \in L^\infty(\R^n \times (0, \infty))$, $V(x,t)= V_1(x,t) +
V_2(x,t)$, $V_2$ is supported in $\{(x,t) : \abs{x} \geq 1\}$,
 and for $C_1,C_2>0$ and $  0\leq \alpha < 1/2$
\begin{equation*}
  \abs{V_1(x,t)} \leq \frac{C_1}{(1+\abs{x}^2)^{\alpha/2}}, 
  \quad
  - (\partial_r V_2 (x,t))^{-} \leq \frac{C_2}{\abs{x}^{2\alpha}}.
\end{equation*}
For $u \in C([0,\infty); L^2(\R^n))$ solution to
\eqref{eq:introduzioneEKPV} there exists a constant $\lambda_0 >0$
such that if
\begin{equation*}
  \sup_{t\geq 0} \int_{\R^n} e^{\lambda_0 \abs{x}^p} \abs{u(x,t)}^2
  \, dx < +\infty, 
  \quad 
  \text{ where } p=(4-2\alpha)/3,
\end{equation*}
then $u$ vanishes. 
Some limit results are also provided  when $\alpha=1/2$.
Moreover, their methods give immediately results in the stationary
case for $\alpha < 1/2$, and  the case $\alpha = 1/2$ is studied
separately, by means of the appropriate Carleman estimate.
For solutions
to the evolution equation  \eqref{eq:introduzioneEKPV},
 properties of  unique
continuation from infinity have been investigated 
in many different papers, also in presence of
electromagnetic perturbations, exploiting a connection with the Hardy
uncertainty principle: we refer to 
\cite{escauriaza2010sharp, escauriaza2016hardy,
  escauriaza2012uniqueness, barcelo2013hardy, cassano2015sharp, cassano2017gaussian}
and references therein. 

In \cite{bourgain2005localization, kenig2006some, kenig2007some}
a quantitative approach to the problem was considered:
 for $u$ solution to \eqref{eq:introduzione1}, with $u$, $V$ bounded and
$u(0)=1$,
it was shown that for big $R \gg 1$ and $C_1,C_2>0$
\begin{equation*}
  M(R):= \inf_{\abs{x_0} = R} \norm{u}_{L^2(B_1(x_0))}
\geq C_1 e^{-C_2 R^{4/3} \log R}.
\end{equation*}
This result was generalized in \cite{davey2014some},
where Davey considered the equation
\begin{equation}\label{eq:introduzione2}
  -\Delta u + W \cdot \nabla u + V u = \lambda u \quad \text{ in }\R^n,
\end{equation}
for $\lambda \in \C$, $\abs{V(x} \lesssim (1+ \abs{x}^2)^{-N/2}, \abs{W(x)} \lesssim
 (1+ \abs{x}^2)^{-P/2}$, $N,P\geq 0$.
For $u$ solution to \eqref{eq:introduzione2} bounded and such that
$u(0) \geq 1$, setting $\beta:= \max(2-2P, (4-2N)/3)$, she showed that
for big $R \gg 1$ and $C_1,C_2,C_3> 0$
\begin{equation*}
  M(R) \geq
  \begin{cases}
    C_1\exp[-C_2 R^{\beta} (\log R)^{C_3}] \quad 
    &\text{ if } \beta >1,
    \\
    C_1\exp[-C_2 R (\log R)^{C_3 \log \log R}] \quad 
    &\text{ if } \beta <1.
  \end{cases}
\end{equation*}
Moreover, explicit examples with the critical decay are provided:
if $V$ and $W$ do not both decay too quickly, these examples are built in the style of
\cite{meshkov1992possible}, otherwise the constructions are simpler.
The case $\beta=1$ was treated in \cite{lin2014quantitative}, where it
is proved that for big $R\gg 1 $ and $C_1,C_2,C_3 >0$ 
\begin{equation*}
  M(R) \geq C_1 \exp[-C_2 R (\log R)^{C_3 (\log R) (\log \log \log  R) (\log \log R)^{-2}}].
\end{equation*}
Finally,  Davey showed in \cite{davey2015meshkov} that this estimate is
sharp, providing an example in the style of \cite{meshkov1992possible}.

The case of a real potential $V$ has been finally addressed in 
\cite{kenig2015landis}, where a quantitative form of Landis' conjecture
is proved in $\R^2$.
 Precisely, for $u$ a real-valued
solution of \eqref{eq:introduzione1} for $E=0$,  $V\geq 0$ 
and $\norm{V}_{L^\infty} \leq 1$, if $u (0) = 1$ 
and $|u (x)| \leq \exp [C_0 \abs{x}]$,
then for a sufficiently large $R\gg 1$
\begin{equation*}
  \inf_{\abs{x_0}=R} 
  \sup_{\abs{x-x_0}\leq 1}
  \abs{u (x)} \geq \exp[-C R \log R],
\end{equation*}
where $C$ depends only on $C_0$.
 Similar estimates for equations with bounded magnetic
potentials are also derived, i.e.~for 
the equations $-\Delta u + W \cdot \nabla u + Vu=0$ 
and $-\Delta u + \nabla (W u) + Vu=0$, 
and the corresponding estimates in exterior domains
$\R^2 \setminus B_R(0)$ are provided.
In \cite{daveykenig2017landis} the results of \cite{kenig2015landis} are generalized 
replacing the Laplace operator with a general operator
$L u:= \text{div}(A \nabla u)$, where 
$A$ is real, symmetric and uniformly elliptic with Lipschitz
continuous coefficients.
Finally, in \cite{davey2017landis} Davey and Wang address the case of
a general second order elliptic equation with singular lower order
terms in $\R^2$.

Regarding the analogous question for the Dirac Operator, to the best of our knowledge, much less results are available.
In \cite{boussaid2016spectral}, Boussa\"id and Comech 
 study the point spectrum of the nonlinear massive Dirac equation
in any spatial dimension, linearized at one of the solitary wave
solutions, and consider the presence of a bounded potential decaying at infinity in 
a weak sense. Thanks to the presence of the massive term, they show
 linear exponential decay
for eigenfunctions and link explicitly the constant in the exponent with the mass, the associated eigenvalue in the
gap of the spectrum
and the ground state for the non-linear stationary equation. 
In the massless case a different behavior should be expected,
as it is shown in \cite{borrelli2018weakly} for the massless non-linear Dirac equation in 
2D without considering the presence of potentials.

In order to state our results we remind the definition and some
elementary properties of the Dirac operator $\D_n$.
For $n\geq 2$ let $N := 2^{\lfloor \frac{n+1}{2}\rfloor}$, where
$\lfloor \cdot \rfloor$ denotes the integer part of a real number.
It is well known that there exist Hermitian matrices
$\alpha_1,\dots,\alpha_n, \alpha_{n+1} \in \C^{N\times N}$
that satisfy the anticommutation relations
\begin{equation}
  \label{eq:anticommutation}
  \alpha_j \alpha_k +
  \alpha_k \alpha_j =
  2 \delta_{j,k} \mathbb{I}_N, 
  \quad 
  1\leq j,k \leq n+1,
\end{equation}
where $\delta_{j,k}$ is the Kronecker delta and $\mathbb{I}_N$ is the 
$N \times N$ identity matrix.
These matrices form a representation of the Clifford algebra of the
Euclidean space (see e.g.~\cite{friedrich2000dirac,jost2008riemannian}):
an explicit choice for them
 is not in general required, since all the results can be obtained
using their fundamental property \eqref{eq:anticommutation};
however it is convenient for our purposes to give one possible choice in the
cases $n=2,3$, that is in the cases where we construct explicitely
examples of solutions to \eqref{eq:main}.
For any $A \in \C^n$ we use the notation
\begin{equation*}
  \alpha \cdot A := \sum_{j=1}^n A_j \alpha_j \in \C^{N \times N}.
\end{equation*}
As an immediate consequence of \eqref{eq:anticommutation} we have that
\begin{equation}
  \label{eq:alpha.A}
  (\alpha \cdot A)^2 = (A \cdot A) \mathbb{I}_N,
  \quad
  \text{ for all }A \in \C^n.
\end{equation}

The Dirac operator acts on functions $\psi : \R^n \to \C^{N}$ and
it is defined as follows:
\begin{equation}
  \D_{n,m} \psi := -i \alpha \cdot \nabla \psi + m \alpha_{n+1}\psi
  := -i \sum_{j=1}^n \alpha_j \partial_j \psi + m \alpha_{n+1}\psi,
\end{equation}
where $m \in \R$. From \eqref{eq:anticommutation} we have that
$\D_{n,m}^2 = (-\Delta + m^2)\mathbb{I}_{N}$.
We remark that  $\D_{n,m}$ is self-adjoint in $L^2(\R^n; \C^N)$
with domain $H^1(\R^n; \C^N)$. 
In the following, when $m=0$ we will denote $\D_{n}:=\D_{n,0}$.

For $n = 2$, we choose explicitely
$(\alpha_1,\alpha_2,\alpha_3):=(\sigma_1,\sigma_2,\sigma_3)$, where
$\sigma_j$ are the \emph{Pauli matrices}
\begin{equation}\label{eq:paulimatrices}
\quad{\sigma}_1 =
\begin{pmatrix}
0 & 1\\
1 & 0
\end{pmatrix},\quad {\sigma}_2=
\begin{pmatrix}
0 & -i\\
i & 0
\end{pmatrix},
\quad{\sigma}_3=
\begin{pmatrix}
1 & 0\\
0 & -1
\end{pmatrix}.
\end{equation}
In the following, when $n=2$, we will sometimes write $\sigma\cdot A:=
\alpha\cdot A$, for any $A\in \C^2$. 
 The two-dimensional Dirac operator is
\begin{equation*}
	\D_{2,m} := 
	- i \sigma_1 \partial_1 - i \sigma_2 \partial_2 + m \sigma_3 =
	\begin{pmatrix}
	m & -2i \partial_z \\
	-2i \partial_{\bar z} & -m 
         \end{pmatrix},
\end{equation*}
where we denote $2\partial_z = \partial_1 -i \partial_2$,
$2\partial_{\zb} = \partial_1 + i \partial_2$, using the usual identification
$(x,y)\in\R^2 \mapsto z \in \C$, with
$z=x_1 + i x_2$, $\zb = x_1 -i x_2$.
 In the following we also denote $r:= \abs{z} = \sqrt{x^2+y^2}$ and
 make use of polar coordinates in $\R^2$.

 For $n=3$, we choose
 \begin{equation*}
   \alpha_j:=\begin{pmatrix}
0& {\sigma}_j\\
{\sigma}_j&0
\end{pmatrix}\quad \text{for}\ j=1,2,3,
\quad
\alpha_4 := \beta :=
\begin{pmatrix}
  \mathbb{I}_2 & 0 \\
  0 & -\mathbb{I}_2
\end{pmatrix}.
 \end{equation*}

We can now state our first result.
 \begin{theorem}\label{thm:negative}
For $n\geq 2$, let $\psi \in H^1_{loc}(\R^n;\C^N)$ be a solution to 
\begin{equation}\label{eq:2d}
  \D_n \psi + \V \psi = 0,
\end{equation}
where $\V : \R^n \to \C^{N \times N}$ is such that
for $-\infty<\epsilon<1$ and $\rho,C>0$
  \begin{equation}\label{eq:askonV}
    \abs{\V(x)} \leq C \abs{x}^{-\epsilon}, \quad \text{ as
    }\abs{x}>\rho.    
  \end{equation}
  Then there exists $\tau_0 >0$ such that, if for all $\tau > \tau_0$
   \begin{equation}\label{eq:conditionu}
     e^{\tau \abs{x}^{2-2\epsilon}} \psi \in L^2(\{\abs{x}>\rho\};\C^N),
   \end{equation} 
then $\psi$ has compact support.
\end{theorem}
\begin{remark}\label{rem:massive}
  In the case that $\V$ verifies the condition \eqref{eq:askonV} for
  $\epsilon \leq 0$, we can also treat solutions
   $\psi \in H^1_{loc}(\R^n;\C^N)$ to 
   \begin{equation}\label{eq:massive.equation}
     \D_{n,m} \psi
     + \V \psi = E \psi,
   \end{equation}
   for $m \in \R$ and $E \in \C$,
   since $\widetilde\V:=\V + m \alpha_{n+1} - E$ 
   satisfy the condition \eqref{eq:askonV} if
  $\epsilon \leq 0$. Thanks to \Cref{thm:negative},
   we conclude that an eigenfunction of the 
    massive Dirac operator perturbed with such a potential $\V$
    can not decay at at infinity more than
    $\exp[-\tau_0 \abs{x}^{2-2\epsilon}]$, in the $L^2$ sense,
    unless it has compact support.
    When $0 < \epsilon < 1$, \Cref{thm:negative} only implies that 
    a solution to \eqref{eq:massive.equation} 
    has compact support if it decays more than a gaussian.
 \end{remark}
 \begin{remark}
   In the assumptions of
   \Cref{thm:negative}
   one would like to conclude that in fact
   the function $\psi$   vanishes everywhere.
   This is possibile assuming some informations on the local behaviour
   of $\V$
   that imply the validity of the \emph{unique continuation} property for
   the operator $\D_n + \V$:
   an operator $A$ is said to satisfy the {unique continuation
     property} if solutions to $A u = 0$ that vanish in a non-empty,
   open subset of a connected set vanish identically in it.
   A vast literature is available on the topic:  we refer to 
   \cite{hile1976unique,kalf1981non,jerison1986carleman,vogelsang1987absence,
     berthier1987point,mandache1994some} and to the survey
   \cite{kenig1989restriction}. We underline here
   two papers among these, that complete the result of \Cref{thm:negative}.
   For $n\geq 2$, in \cite{hile1976unique} it is shown that
   $\D_n + \V$ has the unique continuation property if 
   $\V \in L_{loc}^{\infty}(\R^n;\C^{N \times N})$.
   For $n \geq 3$, in \cite{jerison1986carleman} 
   it is shown that $\V \in L_{loc}^{(3n-2)/2}(\R^n;\C^{N \times N})$
   is enough to guarantee the unique continuation property for the
   operator $\D_n+\V$.
   With these additional assumptions, we can conclude that if $\psi$ has
   compact support, then it vanishes everywhere.
 \end{remark}
 \begin{remark}\label{rem:UC}
   For $\epsilon=1$ the potential $\V$ is a critical perturbation of the Dirac operator:
   we expect that an analogous of \Cref{thm:negative} will employ
   polynomial weights in place of exponential weights at infinity.
   Indeed (see \cite[Remark 1.11]{cassano2018self}), in $\R^3$ and 
   for $\V(x) = -1/\abs{x}$, $m>0$,  the ground state of the 
   Dirac-Coulomb operator is the function
  \begin{equation*}
    \psi_0(x) =    
    \frac{e^{-m\abs{x}}}{\abs{x}}
    \begin{pmatrix}
            1 \\ 1 \\
      i  \sigma \cdot \hx \cdot
      \begin{pmatrix}
        1 \\ 1
      \end{pmatrix}
      \\
    \end{pmatrix},
  \end{equation*}
i.e.~$\psi_0$ is solution to the equation 
   $(-i \alpha\cdot \nabla + m\beta - {\nu}/{\abs{x}} )\psi = 0$.
In the massless case $m=0$, $\psi_0$ is not in $L^2(\R^3)$.
 We refer to \cite{cassano2018self,cassano2019boundary} for a description of the
 functions in the domain of $\D_3 + \V$: we hope to investigate the
 phenomena described in this paper in the case that $\epsilon \geq 1$ in the
 future.
 \end{remark}
 \begin{remark}
   We observe that our results fit in the known theory available for the Laplace operator
   if we ask more regularity on the potential $\V$ and the function $u$.
   For example, let us consider $\V \in W^{1, \infty}(\R^3;\C^{4\times 4})$ and $u \in H^2_{loc}(\R^3;\C^4)$,
   solution to \eqref{eq:2d}: then $u$ is solution to
   \begin{equation*}
     -\Delta u + \sum_{j=1}^3 - i (\alpha_j \cdot  \V)  \partial_j u  -i (\alpha_j \partial_j \cdot \V) u  =0.
   \end{equation*}
   Thanks to \cite{davey2014some}, $u$ has compact support if $\exp[\tau \abs{x}^2] u \in L^2(\R^3)$, for a big 
   enough $\tau \gg 1$, in accordance with \Cref{thm:negative}.
 \end{remark}


The following theorem shows that \Cref{thm:negative} is sharp for $n=2,3$.
\begin{theorem}\label{thm:controesempio}
For all $\epsilon < 1$ and $n\in\{2,3\}$  there exist nontrivial functions
  \begin{equation}\label{eq:propuV2d} 
    u \in C^{\infty}(\R^n;\C^N),  \quad \mathbb{V} \in C^{\infty}(\R^n;\C^{N\times N}),
  \end{equation}
    such that  for all $x \in \R^n$
  \begin{align}
      \label{eq:main2d}
      & \D_n u(x) = \mathbb{V}(x) u(x), \\
    \label{eq:decayu2d}
      & \abs{u(x)} \leq C_1 e^{-C_2 \abs{x}^{2-2\epsilon}}, \\
    \label{eq:decayV2d}
      & \abs{\mathbb{V}(z)}\leq
        \begin{cases}
        C_3 \abs{x}^{-\epsilon}, \quad &\text{ for }n=2,
        \\
        C_3  \abs{x}^{-\epsilon}(\log{\abs{x}})^3,
        \quad &\text{ for }n=3,
        \end{cases}
  \end{align}
  for some $C_1,C_2,C_3>0$.
\end{theorem}
\begin{remark}
  In \cite{duyckaerts2008optimality} the authors 
  construct an explicit function in $\R^3$  with critical decay
  for the Laplace operator and for logaritmically growing potentials,
  moreover they suggest that an analogous construction can be done in $\R^2$:
  their examples are respectively $\C^4$ and $\C^2$--valued.
  The critical examples we construct in 
  \Cref{thm:controesempio} 
  are similar to the ones in
  \cite{duyckaerts2008optimality} and somehow motivate this numerology.
  In detail, in the case $\epsilon=0$ and $n=3$
  the function we build in \Cref{thm:controesempio} is similar to 
  the example given in \cite[Theorem 3.2]{duyckaerts2008optimality};
  the generalization to the case $\epsilon < 1$  is done in light of the 
  approach of \cite[Theorem 2]{cruz1999unique}.
\end{remark}
\begin{remark}
  In accordance with the established theory for the Laplace operator, 
  in \Cref{thm:controesempio}
  the constructed potentials $\V$ are not symmetric and the functions $u$ are complex-valued.
  It would be interesting to understand what is the sharp decay associated to symmetric potentials.
\end{remark}
\subsection*{Structure of the paper}
The paper is organised as follows: 
the proof of \Cref{thm:negative} is given in \Cref{sec:proof-thm:negative};
\Cref{thm:controesempio} is proved in \Cref{sec:controesempio2d}
in the case that $n=2$ and in \Cref{sec:controesempio3d} for $n=3$.

\subsection*{Acknowledgements}
We wish to thank Luca Fanelli and Luis Vega
for addressing me to this theme of investigation, and for precious
discussions and advices.
Moreover, we wish to thank the anonymous referee for noticing a
mistake in the proof of a previous version of \Cref{lem:carleman.astratta}.

\section{Proof of \Cref{thm:negative}}
\label{sec:proof-thm:negative}

In the proof of \Cref{thm:negative} we exploit the following Carleman
Estimate, proved in \cite{enblom2015hardy}.
Carleman estimates have been developed in the study of the Laplace operator,
and there have been many contributions regarding the Dirac operator, mainly 
because of their application to prove unique continuation. We refer to 
\cite{amrein1982estimations, jerison1986carleman, de1990optimal, eller2008carleman,
  salo2009carleman} and references therein for further details.

\begin{proposition}{\cite[Example 4.2]{enblom2015hardy}}
  \label{prop:CarlemanConcreta}
Let $\tau\in \R$ and  $a > 0$.
For all $u\in C^\infty_c(\R^n\setminus\{0\}; \C^N)$ the following holds:
  \begin{equation}\label{eq:carleman}
      \tau a^2 \int_{\R^n} \abs{x}^{a-2} 
      e^{2\tau \abs{x}^{a}} \abs{u(x)}^2 \,dx
    \leq \int_{\R^n} e^{2\tau \abs{x}^{a}} \abs{\D_n u(x)}^2\,dx.
  \end{equation}
\end{proposition}
In fact a Carleman estimate for $\D_n$ is
extablished with general weights in \cite{enblom2015hardy}, and \eqref{eq:carleman} is obtained as a
consequence. For the reader's convenience,
in the following Lemma
we give a short proof of such Carleman estimate
for radial weights, adapting the proof of Theorem 1.1 in \cite{cassano2019hardy}:
we get immediately \eqref{eq:carleman}
setting 
  \begin{equation*}
    b(x):= e^{\tau\abs{x}^{a}}, 
    \quad
    v(x):= b(x) u(x).
  \end{equation*}
\begin{lemma}\label{lem:carleman.astratta}
  Let $b : \R^n \to (0,+\infty)$ be a smooth enough radial function.
  For all $v \in C^\infty_c(\R^n\setminus\{0\}; \C^N)$ we have
  \begin{equation}\label{eq:carlemanAstratta}
    \int_{\R^n} \left(\left(\partial_r^2 +
        \frac{1}{r}\partial_r \right) \log b(x) \right) \abs{v(x)}^2\,dx
    \leq \int_{\R^n} \abs*{b(x) \, \D_n  \, b(x)^{-1} \, v(x)}^2\,dx.
  \end{equation}
\end{lemma}
\begin{remark}
  For $n=2$, \eqref{eq:carlemanAstratta} is a consequence of the analogous Carleman estimate
  for the $\overline\partial$ operator, see \cite[Theorem
  15.1.1]{hormanderII} and \cite[Proposition 2.1]{donnelly1990nodal}:
  the weight that appears at left hand side is
  $\Delta \phi$, when the operator $\overline\partial$ is conjugated
  with the weight $\phi$. Indeed, 
  for $n=2$ and for radial functions $f$, we have $(\partial_r^2 + r^{-1}\partial_r) f =
  \Delta f$.
\end{remark}
\begin{proof}[Proof of \Cref{lem:carleman.astratta}]
  We exploit a radial decomposition of the Dirac operator. We denote $r := \abs{x}$,
  $\widehat{x}:= x/\abs{x}$, and 
  we recall that (see e.g.~\cite{enblom2015hardy})
  \begin{equation}\label{eq:dirac.radial.decomposition.general}
    \D_n = -i \alpha \cdot \nabla
    =
    -i \alpha \cdot \widehat{x}
    \left(\left(
      \partial_r + \frac{n-1}{2 r}\right) \mathbb{I}_N + \frac{S}{r}
    \right),
  \end{equation}
  with 
  \begin{equation*}
    S := \sum_{1\leq j < k \leq n}
    \alpha_j \alpha_k (x_j \partial_k - x_k \partial_j)
    -\frac{n-1}{2}.
  \end{equation*}
  The operator $\partial_r + \frac{n-1}{2r}$ is anti-symmetric on
  $C^\infty_c (\R^n \setminus \{0\};\C^N)$ and acts only on functions
  depending on the radial
  variable, the operator $S$ is
  symmetric on $L^2(\mathbb{S}^{n-1};\C^N)$ acts only on functions
  depending on the angular variable. Moreover, it is useful to observe that
  $\left[\partial_r + \tfrac{n-1}{2r}, S\right] = 0$.

  Thanks to \eqref{eq:alpha.A}, for every $x \in \R^n \setminus\{0\}$ the matrix $-i \alpha \cdot
  \widehat{x}$ is unitary: we have that
    \begin{equation*}
      \begin{split}
      \int_{\R^n} \abs{ b(x) \D_n b(x)^{-1} v(x)}^2\,dx
      = &
      \int_{\R^n} \abs*{ b(r) \left( \partial_r +\frac{n-1}{2r} +
          \frac{S}{r}   \right) b(r)^{-1} v(x)}^2\,dx
      \\
      = &
      \int_{\R^n} \abs*{\left(\partial_r +\frac{n-1}{2r} +
          \frac{S}{r} - \frac{\partial_r b(r)}{b(r)}  \right)  v(x)}^2\,dx
    \end{split}
  \end{equation*}
  Denoting $\abs{x}^{-1/2} v =: \widetilde v \in C^{\infty}_c(\R^n\setminus\{0\};\C^N)$, the right hand side in the
  previous equation can be rewritten as
  \begin{equation*}
    \begin{split}
      \int_{\R^n} &r \abs*{\left(\partial_r +\frac{n-1}{2r} +
          \frac{S}{r} - \frac{\partial_r b(r)}{b(r)} +\frac{1}{2r}
        \right) \widetilde{v}(x)}^2\,dx
      \\
      = &\int_{\R^n} r \abs*{\left(\partial_r +\frac{n-1}{2r}\right)
        \widetilde{v}(x)}^2\,dx
      +
      \int_{\R^n} r \abs*{\left(\frac{S}{r} - \frac{\partial_r b(r)}{b(r)} +\frac{1}{2r}
        \right)
        \widetilde{v}(x)}^2\,dx
      \\
      & + 2 \Re
      \int_{\R^n} 
      \left\langle \left(\partial_r +\frac{n-1}{2r}\right) \widetilde{v}(x), 
      \left( S +\frac{1}{2} - \frac{\partial_r b(r)}{b(r)}r \right)
      \widetilde{v}(x) \right\rangle_{\C^N}\, dx
    \\
    = & \, I + II + III \geq III.
       \end{split}
     \end{equation*}
     Since $\mathcal{A}:=\partial_r +{n-1}/{2r}$ is anti-symmetric and
     $\mathcal{S}:= S +\frac{1}{2} - \frac{\partial_r b(r)}{b(r)}r$ is symmetric, we
     have that (see \cite[Lemma 2.1]{cassano2019hardy}) 
     \begin{equation*}
       III = \int_{\R^n} \langle [\mathcal{S}, \mathcal{A}] \widetilde{v}(x),
       \widetilde{v}(x) \rangle_{\C^N}\, dx
       = \int_{\R^n} \left\langle \frac{1}{r}\left(\partial_r  \frac{\partial_r b(r)}{b(r)}r
         \right){v}(x),
         {v}(x)
       \right\rangle_{\C^N}\, dx.
     \end{equation*}
     With a final explicit computation we get \eqref{eq:carlemanAstratta}.
\end{proof}

We can now prove \Cref{thm:negative}.

Let $a:=2-2\epsilon$. We divide the proof in two steps.

\subsection{Step 1}
We show that, for $\tau >0$ big enough,
for all $u \in C_c^\infty(\{\abs{x}>\rho\};\C^N)$ the
following holds:
\begin{equation}\label{eq:carlemanV}
  \frac{\tau a^2}{4} 
  \int_{\abs{x}>\rho} e^{2\tau \abs{x}^{a}} \abs{x}^{a-2} 
  \abs{u}^2 
  \, dx
  \leq 
  \int_{\abs{x}>\rho} e^{2\tau \abs{x}^{a}} \abs{(\D_n + \V)u}^2 
  \, dx.
\end{equation}
Indeed, thanks to \Cref{prop:CarlemanConcreta}, for $\tau\in \R$,
 $a=2-2\epsilon >  0$, and $u \in C_c^\infty(\{\abs{x}>\rho\};\C^N)$,
\eqref{eq:askonV} and  the elementary inequality $(a+b)^2 \leq 2(a^2 + b^2)$, for
  $a,b\in \R$, we have
 \begin{equation*}
   \begin{split}
        \tau a^2 
   \int_{\abs{x}>\rho} e^{2\tau \abs{x}^{a}} \abs{x}^{a-2} 
  \abs{u}^2 
  \, dx
  \leq &
  2  \int_{\abs{x}>\rho} e^{2\tau \abs{x}^{a}} \abs{(\D_n + \V)u}^2 
  \, dx
  + 2  \int_{\abs{x}>\rho} e^{2\tau \abs{x}^{a}} \abs{\V u}^2 
  \, dx
  \\
  \leq &
2  \int_{\abs{x}>\rho} e^{2\tau \abs{x}^{a}} \abs{(\D_n + \V)u}^2 
   \, dx
   +
   2 C^2
  \int_{\abs{x}>\rho} e^{2\tau \abs{x}^{a}} \abs{x}^{-2\epsilon} 
  \abs{u}^2 
  \, dx.
   \end{split}
 \end{equation*}
Since $u$ has compact support in $\{ \abs{x} > \rho\}$, the second
term at right hand side is finite and can be subtracted from both
sides. Thanks to the arbitrariness in the choice of $\tau \in \R$, we
have that \eqref{eq:carlemanV} holds for $\tau$ big enough.

\subsection{Step 2}
Thanks to an approximation process, \eqref{eq:carlemanV} holds
for all $u \in H^1_{loc}(\{\abs{x}>\rho\};\C^N)$.

Let $h$  be a $C^\infty$ function such that $h(r) = 0$ for $r \leq
1$ and $h(r)= 1$ for $r \geq 2$, and $h_\rho(r):=h(r/\rho)$. Let moreover
$k$ be a $C^\infty_c$ function such that $k(r)=1$ for $r \leq 1$ and
$\text{supp } k \subset \{\abs{x}\leq 2\}$, and $k_R(r):= k (r/R)$.

Set $u_R(x) := h_\rho(\abs{x}) k_R (\abs{x}) \psi(x)$ for all $x \in \R$:
such function is in $H^1(\{\abs{x}>\rho\};\C^N)$ and has compact support.
Since $(\D_n + \V) \psi=0$,
we get, for the appropriate $C>0$, 
\begin{equation*}
\begin{split}
    \frac{\tau a^2}{4} &
  \int_{\abs{x}>\rho} e^{2\tau \abs{x}^{a}} \abs{x}^{a-2} 
  \abs{h_\rho}^2 \abs{k_R}^2 \abs{\psi}^2 
  \, dx
  \leq 
  2 \int_{\abs{x}>\rho} e^{2\tau \abs{x}^{a}} \abs{\D_n(h_\rho k_R)}^2
  \abs{\psi}^2  
  \, dx 
  \\
  &\leq
  \frac{C}{\rho^2}  \int_{\abs{x}>\rho} e^{2\tau \abs{x}^{a}} 
   \chi_{\{\rho<\abs{x}<2\rho\}} 
  \abs{\psi}^2  
  \, dx 
  +
  \frac{C}{R^2}  \int_{\abs{x}>\rho} e^{2\tau \abs{x}^{a}} 
  \chi_{\{R<\abs{x}<2R\}} 
  \abs{\psi}^2  
  \, dx.
\end{split}
\end{equation*}
Letting $R \to +\infty$ and for the appropriate $C,C'>0$, we have that
\begin{equation*}
  \begin{split}
    {\tau a^2} 
    \int_{\abs{x}>\rho} e^{2\tau \abs{x}^{a}} \abs{x}^{a-2} 
    \abs{h_\rho}^2  \abs{\psi}^2 
    \, dx
    \leq &
    \frac{C}{\rho^2}  \int_{\abs{x}>\rho} e^{2\tau \abs{x}^{a}} 
    \chi_{\{\rho<\abs{x}<2\rho\}} 
    \abs{\psi}^2  
    \, dx
    \\
    \leq &
    {C'}  \int_{\rho<\abs{x}<2\rho} e^{2\tau \abs{x}^{a}} 
    \abs{x}^{-2}
    \abs{\psi}^2  
    \, dx.
  \end{split}
\end{equation*}
From the monotonicity of $e^{\tau \abs{x}^a}$, and since
$\abs{x}^{-2} \leq \abs{x}^{a-2}$
\begin{equation*}
  {\tau a^2}
  e^{2\tau \abs{2\rho}^{a}}
  \int_{\abs{x}>2\rho} \abs{x}^{a-2} 
  \abs{\psi}^2 
  \, dx
  \leq
    {C'} e^{2\tau \abs{2\rho}^{a}}   
    \int_{\rho<\abs{x}<2\rho} 
    \abs{x}^{a-2}
    \abs{\psi}^2  
    \, dx.
\end{equation*}
Simplifying the term $ e^{2\tau \abs{2\rho}^{a}}$ in both sides,
the right hand side is smaller than a constant independent of $\tau$:
 from the arbitrariness of $\tau \in \R$ we get that $\psi \equiv 0$
in $\{\abs{x}>2\rho\}$, that is the thesis.

\section{Proof of \Cref{thm:controesempio} for $n=2$}
\label{sec:controesempio2d}
Let $n=2$. The strategy of the proof of \Cref{thm:controesempio} is the
following:
\begin{itemize}
\item we break down $\R^2$ in annuli $\{\rho_k \leq \abs{z}\leq
  \rho_{k+1}\}$, for the appropriate $\rho_k > 0$, $\rho_k \to
  +\infty$;
in such annuli we define the functions $\E_k$ (see
  \Cref{sec:preliminary.defn});
\item for $k\in \N$ big enough, we define the functions $u_k$ and
  $\V_k$ in the annulus $\{\rho_k \leq \abs{z}\leq
  \rho_{k+1}\}$  (see \Cref{sec:defn.annulus});
\item we define the functions $u$ and $\V$ in   $\{\abs{z}\leq
  \rho_{k_0}\}$, for $k_0$ big enough (see \Cref{sec:defn.origin});
\item  in $\{\abs{z}\geq \rho_{k_0}\}$ we define $u$ and $\V$ glueing together the
  functions $u_k$ and $\V_k$, for $k\geq k_0$ (see \Cref{sec:defn.glueing});
\item we check the behaviour at infinity of the function $u$ (see \Cref{sec:decay.u}).
\end{itemize}

\subsection{Notation}
In the following we will write, for sequences of real numbers $(A_k)_{k\in\N}$, $(B_k)_{k\in\N}$,
\begin{equation*}
  A_k = O( B_k)  
\end{equation*}
when there exist constants C>0 and $\bar k$ such that 
\begin{equation*}
  \forall k \geq \bar k, \quad \abs{A_k} \leq C \abs{B_k}.
\end{equation*}
When $A_k$ and $B_k$ also depend on $r \in I$, being $I$ an interval, the estimate is also
assumed to be uniform with respect to $r \in I$.
We will also use the notation
\begin{equation*}
  A_k \approx B_k \iff (A_k=O(B_k) \text{ and }B_k=O(A_k) ).
\end{equation*}

\subsection{Preliminary definitions}\label{sec:preliminary.defn}
In this section we collect some definitions and elementary results we
 need in the proof.

\subsubsection{Definition of $\delta,n_k, d_k, \rho_k, a_k$}
Let 
\begin{equation}
 \label{eq:relazioneepsdelta}
        \delta := 1 - 2\epsilon.
\end{equation}
 Let $n_0$ be a large odd number, and for all $k \geq 0$ define the following quantities: 
\begin{align}
  \label{eq:propnk}
  &n_{k+1} := n_{k} + 2 \lfloor n_k^{1/2}\rfloor, 
  \quad d_k:= \frac{n_{k+1} - n_{k}}{2},
  \\
  \label{eq:proprhok}
  &\rho_{k}:= n_k^{\frac{1}{1+\delta}},
  \quad
  \rho_{kj}:=\rho_{k} + j \, \frac{\rho_{k+1}-\rho_{k}}{4}, \quad j=1,\dots,4,
  \\
  \label{eq:defnak}
  &a_0:= 1, 
  \quad
  a_{k+1}:= \rho_k^{2d_k} a_k.
\end{align}

It is easy to see that 
\begin{gather}
  \label{eq:propdk}
   d_k = n_k^{1/2}+ O(1),
   \\
   \label{eq:differenzarhok}
  \rho_{k+1}- \rho_{k} 
  = 
  n_{k+1}^{\frac{1}{1+\delta}} - n_{k}^{\frac{1}{1+\delta}} 
   =  
  \frac{1}{1+\delta} n_k^{-\frac{\delta}{1+\delta}} (2 n_k^{1/2}+ O(1))
  = \frac{2}{1+\delta} \rho_k^{\frac{1-\delta}{2}} + O(\rho_k^{-\delta}).
\end{gather}
Moreover, we show that
\begin{equation}\label{eq:at.infinity}
  k = O(\rho_k^{\frac{1+\delta}{2}}),
\end{equation}
i.e.~there exist $k_0 \in \N$ and $C>0$ such that
$  k \leq C \rho_k^{\frac{1+\delta}{2}}$
for all $k \geq k_0$.
We show this by induction:
setting $C :=  \text{max}(k_0 \rho_{k_0}^{- \frac{1+\delta}{2}},3)$, the
basis of induction is true.
Assume now that the thesis is true for some $k \geq k_0$.
 From \eqref{eq:propnk} and \eqref{eq:proprhok}, assuming that $k_0$ is big enough, we have that 
\begin{equation*}
  \rho_{k+1}^{1+\delta} = \rho_k^{1+\delta}
  + 2 \lfloor \rho_k^{\frac{1+\delta}{2}} \rfloor
  \geq \rho_k^{1+\delta} + \rho_k^{\frac{1+\delta}{2}}.
\end{equation*}
We conclude immediately that the thesis is valid for $k+1$: 
\begin{equation*}
  C \rho_{k+1}^{\frac{1+\delta}{2}}
  \geq
  (C^2 \rho_{k}^{1+\delta}
  +
  C^2 \rho_{k}^{\frac{1+\delta}{2}})^{\frac12}
  \geq
  (k^2 + C k)^{\frac12}
  \geq (k^2 + 2 k + 1)^{\frac12} =  k+1,
\end{equation*}
the last inequality being true for $k_0$ big enough.

Finally,  for  $\rho_k \leq r \leq \rho_{k+1}$ we have that
\begin{equation}\label{eq:equivrk}
  \frac{a_k}{r^{n_k}} \approx \frac{a_{k+1}}{r^{n_{k+1}}}.
\end{equation}
Indeed, denoted $g(r):=\log(r^{d_k}/\rho_k^{d_k})$ for 
$\rho_k \leq r \leq \rho_{k+1}$, thanks to
\eqref{eq:propdk} and \eqref{eq:proprhok} we have that 
\begin{equation*}
  g'(r) = \frac{d_k}{r} = \frac{O(\rho_k^{\frac{1+\delta}{2}})}{r}
        = O(\rho_k^{-\frac{1-\delta}{2}}).
\end{equation*}
Observing that $g(\rho_k)=0$ and $\rho_{k+1}-\rho_{k} = O(\rho_k^{\frac{1-\delta}{2}})$,
we get that $g(r)=O(1)$,
 that is
\begin{equation}\label{eq:last}
  r^{d_k}\approx \rho_k^{d_k}, \quad \text{ for }\rho_k \leq r \leq \rho_{k+1}.
\end{equation}
From \eqref{eq:last} we get immediately \eqref{eq:equivrk}.

\subsubsection{Definition of $\E_k$}
For $k\in\N$ big enough and $r:=\abs{z}\geq \rho_k$ we define
\begin{gather*}
  E_k(z):= \frac{a_k}{z^{n_k}},
  \\
  \E_k(z) :=
  \begin{pmatrix}
    E_k(z) \\ 0
  \end{pmatrix}
  \text{ if $k$ is even}, 
  \quad
  \E_k(z) :=
  \begin{pmatrix}
    0 \\ \overline{E_k(z)}
  \end{pmatrix}
  \text{ if $k$ is odd}. 
\end{gather*}
We observe that for all $ z \in \R^2\setminus\{ 0 \}$ and $k \in \N$:
\begin{equation}\label{eq:zeroEk}
  \D_2 \mathbf{E}_k (z) = 0.
\end{equation}
We remark that \eqref{eq:defnak} implies
$\abs{E_k(z)} = \abs{E_{k+1}(z)}$  whenever $r=\abs{z} =
\rho_{k+1}$. 
Moreover, since $\abs{E_k(z)}=a_k r^{-n_k}$ we have immediately from \eqref{eq:equivrk} that
\begin{equation}\label{eq:equivEk}
   E_{k}(z) \approx E_{k+1}(z), \quad \text{ for } \rho_k \leq \abs{z} \leq \rho_{k+1}.
\end{equation}

\subsection{Definition of $u_k$ and $\V_k$ in the annulus $\{\rho_k \leq
  \abs{z}\leq\rho_{k+1}\}$}
\label{sec:defn.annulus}
For $k\in\N$ big enough, in this section we construct functions 
\begin{equation}\label{eq:regolarityannulus}
  \begin{split}
    &u_k \in C^{\infty}(\{z \in \R^2 \colon \rho_k \leq \abs{z} \leq
    \rho_{k+1} \}; \C^2), \\
    &V_k \in C^\infty(\{z \in \R^2 \colon \rho_k \leq \abs{z} \leq
    \rho_{k+1} \}; \C^{2\times 2}),
  \end{split}
\end{equation}
such that 
\begin{equation}
  \label{eq:main2dlocal}
  \mathcal D_2 u_k (z) = \V (z) u_k(z), \quad \text{ for all }\rho_k
  \leq \abs{z} \leq \rho_{k+1},
\end{equation}
and
\begin{equation}
  \label{eq:ugualeinintervalliestremi}
  u_k(z)= 
  \begin{cases}
    {\mathbf{E}_{k}(z)} \quad &\text{for $r \in [\rho_{k0},\rho_{k1}]$}, \\
    {\mathbf{E}_{k+1}(z)} \quad &\text{for $r \in [\rho_{k3},\rho_{k4}]$},
  \end{cases}
\end{equation}
 for $\rho_k \leq \abs{z} \leq \rho_{k+1}$
\begin{equation}
  \label{eq:decaylemma}
  \abs{u_k(z)} = O (a_k r^{-n_k}),
\end{equation}
and  \eqref{eq:decayV2d} holds.
In the proof in this section we assume $k$ to be a odd integer: the
case of even $k$ can be treated analogously.

\subsubsection{Construction of the cut-off functions}
\label{sec:defn.cut-off}
Let $\chi$ be a $C^\infty$ non-decreasing function on $\R$ such that
$\chi(s)=0$ for $s\leq 0$, $\chi(s)=1$ for $s \geq 1$. For 
$\rho_k \leq r \leq \rho_{k+1}$ let
\begin{equation}\label{eq:defnchik}
  \chi_k(r):= \chi\left( \frac{4}{\rho_{k+1}-\rho_k}  (r-\rho_{k1})  \right), 
  \quad
  \tilde \chi_k(r):= \chi\left(  \frac{4}{\rho_{k+1}-\rho_k}  (\rho_{k3}-r)  \right).
\end{equation}
Thanks to \eqref{eq:differenzarhok},  we have that
\begin{equation}
\label{eq:controlloderivate}
  \norm{\chi_k'}_{L^\infty([\rho_k,\rho_{k+1}])}
  =
  \norm{\tilde \chi_k'}_{L^\infty([\rho_k,\rho_{k+1}])}
  =
  O(\rho_k^{-\frac{1-\delta}{2}}).
\end{equation}

\subsubsection{Definition of the functions $u_k$ and $\mathbb V_k$}
For $\rho_k \leq \abs{z} \leq \rho_{k+1}$, let
\begin{equation*}
\begin{split}
 & u_k(z):= 
  \tilde \chi_k(r) \mathbf{E}_{k}(z)  +
  \chi_k(r) \mathbf{E}_{k+1}(z)
  =
  \begin{pmatrix}
    \chi_k(r) E_{k+1}(z) \\
    \tilde \chi_k(r) \overline{E_{k}(z)}
  \end{pmatrix},
\\
 & \V_k (z) := 
  \begin{cases}
    \begin{pmatrix}
      0 & 0 \\
      0 & 0
    \end{pmatrix}
    \quad &\text{ for }r\in [\rho_{k0},\rho_{k1}]\cup
    [\rho_{k3},\rho_{k4}],
    \\
    \begin{pmatrix}
      0 & 0 \\
      0 & -2i \partial_\zb (\chi_k(r))
{E_{k+1}(z)}({\overline{E_{k}(z)}})^{-1} 
    \end{pmatrix}
    \quad &\text{ for } r\in [\rho_{k1},\rho_{k2}],
    \\
    \begin{pmatrix}
      -2i \partial_z(\tilde\chi(r)) \overline{E_k(z)} (E_{k+1}(z))^{-1} & 0 \\
      0 & 0
    \end{pmatrix}
    \quad &\text{ for }r\in [\rho_{k2},\rho_{k3}].
  \end{cases}
\end{split}
\end{equation*}
By construction \eqref{eq:regolarityannulus}, \eqref{eq:main2dlocal}
and \eqref{eq:ugualeinintervalliestremi} are true;
thanks to \eqref{eq:equivEk}, \eqref{eq:decaylemma} holds 
for a large enough $k$.
Some more details are in order for checking condition 
 \eqref{eq:decayV2d}. 

If $r \in [\rho_{k0},\rho_{k1}]\cup
    [\rho_{k3},\rho_{k4}]$ 
 \eqref{eq:decayV2d} is trivially verified.
In the case $r \in [\rho_{k1},\rho_{k2}]$,
thanks to 
\eqref{eq:equivEk},
\eqref{eq:controlloderivate} and \eqref{eq:relazioneepsdelta},
\begin{equation*}
  \abs{\mathbb{V}_{k}(z)} 
  =O(\rho_k^{-\frac{1-\delta}{2}})
  = O(\rho_k^{-\epsilon}).
\end{equation*}
Thanks to \eqref{eq:differenzarhok}, we get that
\begin{equation}\label{eq:canbeofhelp}
  \frac{\rho_{k+1}}{\rho_k}= O(1 + 4 \rho_k^{-\frac{1+\delta}{2}}) = O(1),
\end{equation}
and we conclude \eqref{eq:decayV2d} for $r \in [\rho_{k1},\rho_{k2}]$  since
\begin{equation*}
    \abs{\mathbb{V}_{k}(z)} 
    = O (\rho_{k+1}^{-\epsilon}) 
    = O (\abs{z}^{-\epsilon}).
\end{equation*}
Finally \eqref{eq:decayV2d} is proved analogously in the interval $[\rho_{k2},\rho_{k3}]$.

\subsection{Definition of $u$ in $\{ \abs{z} \leq \rho_{k_0} \}$}
\label{sec:defn.origin} 
Let $k_0$ be a large integer such that the arguments in the previous section hold
for all $k \geq k_0$.
Without loss of generality we can assume $k_0$
to be odd. 
For $\abs{z} \leq \rho_{k_0}$ let 
\begin{equation*}
  \psi(z):= \chi \left(\frac{4}{\rho_{k_0}}\left(\abs{z} -
      \frac{\rho_{k_0}}{4}\right)\right),
  \quad
  \tilde \psi(z):= \chi \left(\frac{4}{\rho_{k_0}}\left(
      \frac{3\rho_{k_0}}{4}-\abs{z} \right)\right),
\end{equation*}
with $\chi$ defined in \Cref{sec:defn.cut-off} and
\begin{equation}
  \label{eq:1}
  u(z) :=  \tilde \psi(z)
  \begin{pmatrix}
    a_{k_0} z^{n_{k_0}} 
    \\
    0
  \end{pmatrix}
  +
  \psi(z)  \E_{k_0}
  =
  \begin{pmatrix}
    \tilde \psi(z) a_{k_0} z^{n_{k_0}} 
    \\
    \psi(z) a_{k_0}\bar z^{-n_{k_0}}
  \end{pmatrix}.
\end{equation}
We remark that $\D u(z) = 0 $ for $\abs{z}\leq
\rho_{k_0}/4$. Moreover, for 
\begin{equation}
  \V (z) := 
  \begin{cases}
    \begin{pmatrix}
      0 & 0 \\
      0 & 0
    \end{pmatrix}
    \quad &\text{ for }\abs{z}\in [0,\rho_{k_0}/4]\cup
    [3\rho_{k_0}/4,\rho_{k_0}],
    \\
    \begin{pmatrix}
      -2i \partial_z(\psi(z)) \abs{z}^{-2n_{k_0}} & 0 \\
      0 & 0
    \end{pmatrix}
    \quad &\text{ for } \abs{z}\in [\rho_{k_0}/4,2\rho_{k_0}/4],
    \\
    \begin{pmatrix}
      0 & 0 \\
      0 & -2i \partial_\zb (\tilde \psi(z))
      \abs{z}^{2n_{k_0}}
    \end{pmatrix}
    \quad &\text{ for }\abs{z}\in [2\rho_{k_0}/4,3\rho_{k_0}/4],
  \end{cases}
\end{equation}
conditions \eqref{eq:propuV2d} -- \eqref{eq:decayV2d} hold for the
appropriate $C_1, C_2, C_3 >0$.

\subsection{Definition of $u$ in $\{ \abs{z} \geq \rho_{k_0} \}$}
\label{sec:defn.glueing} 
For any $z \in \R^2$, $\abs{z} \geq \rho_{k_0}$, let  $k$ be such that
$\abs{z} \in [\rho_k, \rho_{k+1}]$. We set
\begin{equation}
  u(z) := u_k(z), \quad 
  \V (z) := \V_k(z), \quad \text{ for }   \abs{z} \in [\rho_k,
  \rho_{k+1}].
\end{equation}
It is easy to show that conditions 
\eqref{eq:propuV2d}, \eqref{eq:main2d}, \eqref{eq:relazioneepsdelta}, \eqref{eq:decayV2d}
hold, for the appropriate $C_3 >0$.

\subsection{Decay of $u$ at infinity}
\label{sec:decay.u}
We show in detail that \eqref{eq:decayu2d} holds. 
Let $j\in \N$ big enough and $z   \in \R^2$ such that $\rho_j \leq r = \abs{z} \leq \rho_{j+1}$.
Denoting $\frac{r}{\rho_j}=1+h$, thanks to \eqref{eq:proprhok} and \eqref{eq:differenzarhok} we have that $h=O(\rho_j^{-\frac{1+\delta}{2}})$ and 
$n_j h^2 = O(1)$. Using \eqref{eq:decaylemma} we have that
\begin{equation*}
  \begin{split}
  \log \abs{u(z)} - \log \abs{u(\rho_j)} 
  \leq &
    - n_j \log{r} + n_j \log{\rho_j} + O(1)
  =
    -n_j \log{\frac{r}{\rho_j}} + O(1) 
    \\
  \leq &
    -n_j h + \frac{n_j h^2}{2} + O(1)
    =
    -n_j h + O(1).
  \end{split}
\end{equation*}
Moreover, for $m(r):= e^{-\frac{r^{1+\delta}}{1+\delta}}$, we have
\begin{equation*}
  \log{m(r)}-\log{m(\rho_j)} = -\frac{r^{1+\delta}}{1+\delta} + \frac{\rho_j^{1+\delta}}{1+\delta}
  =
  -\frac{\rho_j^{1+\delta}}{1+\delta} \left((1+h)^{1+\delta}-1\right) 
  =
  -\rho_j^{1+\delta} h + O(1).
\end{equation*}
From the previous reasoning, we have immediately that
\begin{equation}\label{eq:asintotica1}
  \log\abs{u(z)} - \log{\abs{u(\rho_j)}} \leq \log{m(r)} - \log{m(\rho_j)} + O(1).
\end{equation}
From \eqref{eq:asintotica1} we have moreover that for all $l \geq
l_0$,  for $l_0$ big enough,
\begin{equation}\label{eq:asintotica2}
  \log\abs{u(\rho_l)} - \log{\abs{u(\rho_{l-1})}} \leq \log{m(\rho_l)} - \log{m(\rho_{l-1})} + O(1).
\end{equation}

Thanks to \eqref{eq:asintotica1} and \eqref{eq:asintotica2}, we have that
\begin{equation*}
  \begin{split}
    \log\abs{u(z)} =
    &
    \log\abs{u(z)} - \log\abs{u(\rho_j)}
    +
    \sum_{l=l_0 +1}^{j}
    \big(\log\abs{u(\rho_l)}-\log\abs{u(\rho_{l-1})}\big)
    +
    \log\abs{u(\rho_{l_0})}
    \\
    \leq & \log{m(r)} + O(j).  
  \end{split}
\end{equation*}
Thanks to \eqref{eq:at.infinity},
taking the exponential at the two sides in the previous equation
we get that, for $l$ big enough and for $C>0$,
\begin{equation*}
  \abs{u(z)} \leq \text{exp}
  \left[ -\frac{r^{1+\delta}}{1+\delta} + C \rho_j^{\frac{1+\delta}{2}}  \right]
\end{equation*}
Consequently, for all $z \in \R^2$ such that $\abs{z} \geq
\rho_{l_0}$, we have
\begin{equation*}
\abs{u(z)}
   \leq \text{exp}
  \left[ -\frac{r^{1+\delta}}{1+\delta} + C r^{\frac{1+\delta}{2}}  \right].
\end{equation*}
We conclude immediately \eqref{eq:decayu2d}, for the appropriate
$C_1,C_2 > 0$,
asking that $l_0$ is big enough.

\section{Proof of \Cref{thm:controesempio} for $n=3$}
\label{sec:controesempio3d}
The proof of \Cref{thm:controesempio} for $n=3$ is very
similar to the one of the case $n=2$,
so we will just sketch it here underlining the main differences.

\subsection{Preliminary definitions}
In this section we collect some definitions and elementary results we need in the proof.
We set 
\begin{equation*}
r=|x|,\quad\hat x = x / |x|\quad\text{and}\quad
L=-ix \wedge \nabla\quad\text{for }x \in \Rt\setminus\{0\}.
\end{equation*}
In the following we will often use polar coordinates.

For $l=0,1,\dots$ and $m=-l,-l+1,\dots,l$,
let $Y_l^m$ be the \emph{spherical harmonics}
\begin{equation*}
  \begin{split}
    &Y_{l}^{m} (\theta ,\varphi )=\sqrt {\frac{2 l +1}{4\pi}  \frac{(l-m)! }{(l +m)!}} 
    \,e^{im\varphi }\, P_{l}^{m}(\cos {\theta }),   \\
    &  Y_{l}^{-m} = (-1)^{m} \overline{Y_l^m}, 
  \end{split}
\end{equation*}
where $P_l^m$ are the \emph{associated Legendre polynomials}
\begin{equation*}
  P_{l}^{m}(x)={\frac {(-1)^{m}}{2^{l} l !}}(1-x^{2})^{m/2}\ {\frac
    {d^{l +m}}{dx^{l +m}}}(x^{2}-1)^{l }.
\end{equation*}
For $\kappa= \pm 1, \pm 2, \dots$,  $ j = \abs{k}  -1/2 = 1/2, 3/2,
\dots$, and $m_j = -j, -j+1, \dots, j$, let $\Y_{\kappa,m_j} $ be the
\emph{spinor harmonics} 
\begin{equation*}
  \Y_{\kappa, m_j} = 
  \begin{cases}
  \frac{1}{\sqrt{2 \kappa - 1}}
  \begin{pmatrix}
    \sqrt{\kappa - \frac12 + m_j} Y_{\kappa -1}^{m_j -1/2} \\
    \sqrt{\kappa - \frac12 - m_j} Y_{\kappa -1}^{m_j +1/2} 
  \end{pmatrix}, \quad  & \kappa \geq 1, \\
  \frac{1}{\sqrt{1 - 2 \kappa}}
  \begin{pmatrix}
    -\sqrt{\frac12 - \kappa - m_j} Y_{-\kappa}^{m_j -1/2} \\
    \sqrt{\frac12 - \kappa + m_j} Y_{-\kappa}^{m_j +1/2} 
  \end{pmatrix}, \quad  & \kappa \leq -1 .
\end{cases}
\end{equation*}
It is well known (see \cite[Section 3.9.4]{thaller2005advanced} or
\cite[Section 4.6.4]{thaller1992dirac}) that
\begin{equation}\label{eq:propautofunzioni}
(\mathbb{I}_2+\sigma\cdot L) \, \Y_{\kappa,m_j}= \kappa \, \Y_{\kappa,m_j},
\end{equation}
where $\sigma=(\sigma_1,\sigma_2,\sigma_3)$ is the vector of the 
\textit{Pauli matrices} (defined in \eqref{eq:paulimatrices}).

\subsubsection{Definition of $\delta,n_k, d_k, \rho_k, a_k$}
Let $\delta>-1$ as in \eqref{eq:relazioneepsdelta}.
Let $n_0$ be a large odd number and for all $k\geq 0$ let $n_k, d_k,
\rho_k$ and $a_k$ be defined as in  \eqref{eq:propnk} --  \eqref{eq:defnak}. 
\subsubsection{Definition of $F_{n_k}$}
For $m \in \N$ set
\begin{equation}\label{eq:defnFnk}
F_{m}(\theta,\phi):= c_{m}  \Y_{-(m-1),\frac12}(\theta,\phi)
=\frac{c_{m}}{\sqrt{4\pi}} 
  \begin{pmatrix}
    -\sqrt{m-1} P_{m-1}^0 (\cos \theta) \\
    \frac{1}{\sqrt{m-1}} e^{i \phi} P_{m-1}^1(\cos \theta)
  \end{pmatrix},
\end{equation}
with $c_{m}>0$ defined later (in \eqref{eq:choiceck}). 
We show that for all large odd $m$ the functions $F_m$ satisfy
\begin{equation}\label{eq:tesiFnk}
  {C}{m^{-\frac74}} \leq \abs{F_m(\theta,\phi)} \leq 1,
\end{equation}
for some $C>0$. Indeed, from \eqref{eq:defnFnk} we have that
\begin{equation*}
  \abs{F_m(\theta,\phi)}^2 = \frac{c_m^2}{4\pi} 
  \left((m-1)\abs{P_{m-1}^0(\cos \theta)}^2 +
    \frac1{m-1}\abs{P_{m-1}^1(\cos \theta)}^2\right).
\end{equation*}
We remind here the following property of the associated Legendre polynomials, proved in \cite[Lemma 3.8]{duyckaerts2008optimality}:
 there exists $ C >0$ such that for all large $l \in 2\N$
\begin{equation}
  \label{eq:propPml}
  \frac{1}{ C l^{3/2}} \leq l(l+1) \abs{P_l^0(x)}^2 +
  \abs{P_l^1(x)}^2
  \leq  C l(l+1),
  \quad \text{ for all } x \in [-1,1].
\end{equation}
Choosing $l=m-1$, with an easy computation we get
\begin{equation*}
  \frac{1}{ C(m-1)^{3/2}m} \leq
  \frac{4\pi}{c_m^2}\abs{F_m(\theta,\phi)}^2 - \frac{1}{m(m-1)}
  \abs{P_{m-1}^1(\cos \theta)}^2 \leq  C (m-1).
\end{equation*}
From \eqref{eq:propPml} and the last equation we have
\begin{equation}\label{eq:useropertesi}
  \frac{1}{ C(m-1)^{3/2}m} \leq \frac{4\pi}{c_m^2} \abs{F_m(\theta,\phi)}^2 \leq  C (m-1) +  \frac{1}{m(m-1)} \abs{P_{m-1}^1(\theta)}^2 \leq  C' m,
\end{equation}
for an appropriate $C'>0$ and big enough $m$.
We set
\begin{equation}
  \label{eq:choiceck}
  c_m:=\sqrt{\frac{4\pi}{ C' m}},
\end{equation}
and we conclude \eqref{eq:tesiFnk} for the appropriate $C>0$, thanks to \eqref{eq:useropertesi}.
\subsubsection{Definition of $\E_{k}$}
For $k\in\N$ big enough and $r:=\abs{x}\geq \rho_k$ we define
\begin{gather*}
  E_{k}(x) :=  a_k r^{-n_k} F_{n_k}(\theta,\phi), \\
  \E_k(x) :=
  \begin{pmatrix}
    E_k(x) \\ 0
  \end{pmatrix}
  \text{ if $k$ is even}, 
  \quad
  \E_k(x) :=
  \begin{pmatrix}
    0 \\ E_k(x)
  \end{pmatrix}
  \text{ if $k$ is odd}. 
\end{gather*}

Thanks to \eqref{eq:dirac.radial.decomposition.general} for $n=3$
 (see also
\cite[eq.~(4.104)]{thaller1992dirac})  and \eqref{eq:propautofunzioni}
we have that
for all $ x \in \R^3\setminus\{ 0 \}$ and $k \in \N$
\begin{equation*}
  - i \sigma \cdot \nabla \, E_k(x) =
  - i \sigma \cdot \hx \,
  \left( \partial_r +\frac{1}{r} - \frac{1+\sigma\cdot L}{r} \right)
  \, E_k(x) = 
  0,
\end{equation*}
that gives immediately
\begin{equation*}
-i \alpha \cdot \nabla \, \E_k(x) = 0.  
\end{equation*}
Moreover, thanks to \eqref{eq:equivrk}, \eqref{eq:canbeofhelp} and \eqref{eq:tesiFnk} we have that
for $k$ big enough 
\begin{equation}
  \label{eq:equivEk3d}
  \frac{\abs{\E_{k+1}(x)}}{\abs{\E_{k}(x)}} = O( n_k^{\frac74}), 
  \quad
  \frac{\abs{\E_{k}(x)}}{\abs{\E_{k+1}(x)}} = O( n_k^{\frac74}).
\end{equation}
for the appropriate $C>0$. 

\subsection{Definition of $u_k$ and $\V_k$ in the annulus $\{\rho_k \leq
  \abs{z}\leq\rho_{k+1}\}$}
For $k\in\N$ big enough, in this section we construct functions 
\begin{equation}\label{eq:regolarityannulus3d}
  \begin{split}
    &u_k \in C^{\infty}(\{x \in \R^3 \colon \rho_k \leq \abs{x} \leq
    \rho_{k+1} \}; \C^4), \\
    &\V_k \in C^\infty(\{x \in \R^3 \colon \rho_k \leq \abs{x} \leq
    \rho_{k+1} \}; \C^{4\times 4}),
  \end{split}
\end{equation}
such that 
\begin{equation}
  \label{eq:main3dlocal}
  \mathcal D_3 u_k (z) = \V_k (z) u_k(z), \quad \text{ for all }\rho_k
  \leq \abs{x} \leq \rho_{k+1},
\end{equation}
and
\begin{equation}
  \label{eq:ugualeinintervalliestremi3d}
  u_k(x)= 
  \begin{cases}
    {\mathbf{E}_{k}(x)} \quad &\text{for $r \in [\rho_{k0},\rho_{k1}]$}, \\
    {\mathbf{E}_{k+1}(x)} \quad &\text{for $r \in [\rho_{k3},\rho_{k4}]$},
  \end{cases}
\end{equation}
 for $\rho_k \leq \abs{x} \leq \rho_{k+1}$
\begin{equation}
  \label{eq:decaylemma3d}
  \abs{u_k(x)} = O (a_k r^{-n_k}),
\end{equation}
and  \eqref{eq:decayV2d} holds.
In the proof in this section we assume $k$ to be a odd integer: the
case of even $k$ can be treated analogously.

\subsubsection{Construction of the cut-off functions}
Let $\chi$ be a real non-decreasing $C^{\infty}$ 
function on $\R$ such that
\begin{equation*}
  \chi(s) = 
  \begin{cases}
    0 \quad & s \leq 0, \\
    e^{-1/s} \quad & 0 \leq s \leq 1/2, \\
    1 \quad & s\geq 1,
  \end{cases}
\end{equation*}
and for $\rho_k \leq r \leq \rho_{k+1}$ 
let  $\chi_k$ and $\tilde \chi_k$ be defined as in
\eqref{eq:defnchik}. 
We remark that \eqref{eq:controlloderivate} holds.
\subsubsection{Definition of the functions $u_k$ and $\mathbb V_k$}
\label{sec:remarkable}
For $\rho_k \leq r=\abs{x} \leq \rho_{k+1}$, let
\begin{equation*}
  u_k(x):= \tilde \chi_k (r) \E_k(x) + \chi_k(r) \E_{k+1}(x) 
  = 
 \begin{pmatrix}
    \chi_k(r) E_{k+1}(x) \\
    \tilde \chi_k(r) {E_{k}(x)}
  \end{pmatrix}.
\end{equation*}
Let
\begin{equation*}
 \V_k (x) := 
    \begin{pmatrix}
      0 & 0 \\
      0 & 0
    \end{pmatrix},
    \quad r\in [\rho_{k0},\rho_{k1}]\cup [\rho_{k3},\rho_{k4}],
\end{equation*}
and for $r\in [\rho_{k1},\rho_{k2}]\cup [\rho_{k2},\rho_{k3}]$ let
$\V_k(x):= \D_3 u_k(x) \, \overline{u_k(x)}^t / \abs{u_k(x)}^2$, that is
\begin{equation*}
 \V_k (x) := 
  \begin{cases}
    \frac{1}{\abs{\chi_k E_{k+1}}^2+\abs{E_k}^2}
    \begin{pmatrix}
      0 & 0 \\
      0 & -i \sigma \cdot \hx 
    \end{pmatrix}
    \begin{pmatrix}
      0 & 0 \\
      \chi_k' E_{k+1} \cdot \chi_k \overline{E_{k+1}}^t
      & \chi_k' E_{k+1} \cdot \overline{E_{k}}^t 
    \end{pmatrix}
    \, &\text{ for } r\in [\rho_{k1},\rho_{k2}],
    \\
    \frac{1}{\abs{\chi_k E_{k+1}}^2+\abs{E_k}^2}
\begin{pmatrix}
      -i \sigma \cdot \hx  & 0 \\
      0 & 0
    \end{pmatrix}
    \begin{pmatrix}
      \tilde \chi_k' E_{k} \cdot \overline{E_{k+1}}^t 
      & \tilde\chi_k' E_{k} \cdot \tilde\chi_k \overline{E_{k}}^t
      \\
      0 & 0
    \end{pmatrix}
    \, &\text{ for }r\in [\rho_{k2},\rho_{k3}].
  \end{cases}
\end{equation*}

By construction \eqref{eq:regolarityannulus3d}, \eqref{eq:main3dlocal}
and \eqref{eq:ugualeinintervalliestremi3d} are true;
thanks to \eqref{eq:equivEk3d}, \eqref{eq:decaylemma3d} holds 
for a large enough $k$.

Condition \eqref{eq:decayV2d} holds obviously for $r \in [\rho_{k0},\rho_{k1}] \cup [\rho_{k3},\rho_{k4}]$.
For the case $r \in [\rho_{k1},\rho_{k2}]$ we need to distinguish various cases:
let $s:= 4(r-\rho_{k1})/(\rho_{k+1}-\rho_k)$ and let us consider first the
case $s \in (0,(\log \rho_k)^{-3/2})$. 
We have that
\begin{equation}\label{eq:V12prima}
\begin{split}
  \abs{\V_{21}(x)} \leq & \frac{\abs{\chi_k'(r)} \abs{E_{k+1}(x)} \abs{\chi_k(r)} \abs{E_{k+1}(x)}}%
  {\abs{\chi_k(r) E_{k+1}}^2 + \abs{E_{k}(x)}^2} 
  \\
  \leq &
  \frac{\abs{\chi_k'(r)} \abs{E_{k+1}(x)} \abs{\chi_k(r)} \abs{E_{k+1}(x)}}%
  {2\abs{\chi_k(r)}\abs{E_{k+1}(x)}\abs{E_{k}(x)}}  
  =
  \frac{\abs{\chi_k'(r)}\abs{E_{k+1}(x)}}{2\abs{E_{k}(x)}}
\end{split}
\end{equation}
and
\begin{equation}\label{eq:V22prima}
\begin{split}
  \abs{\V_{22}(x)} \leq & \frac{\abs{\chi_k'(r)} \abs{E_{k+1}(x)} \abs{E_{k}(x)}}%
  {\abs{\chi_k(r) E_{k+1}}^2 + \abs{E_{k}(x)}^2} 
  \\
  \leq & \frac{\abs{\chi_k'(r)} \abs{E_{k+1}(x)} \abs{E_{k}(x)}}{\abs{E_{k}(x)}^2}
  =
  \frac{\abs{\chi_k'(r)}\abs{E_{k+1}(x)}}{\abs{E_{k}(x)}}.
\end{split}
\end{equation}
We observe that
\begin{equation}\label{eq:stimaderivata}
  \chi'_k(r) =  \frac{4}{\rho_{k+1}-\rho_{k}} \chi'(s) =
  O(\rho_{k}^{-\frac{1-\delta}{2}}) \chi'(s).
\end{equation}
In this interval $\chi$ has an explicit expression: since $\chi'$ is
increasing we get for large $k$ that
 \begin{equation}\label{eq:stimaderivataprima}
   \chi'_k(r) \leq  
   O(\rho_{k}^{-\frac{1-\delta}{2}}) \chi'((\log \rho_k)^{-3/2})
   = O(\rho_{k}^{-\frac{1-\delta}{2}}) 
   (\log\rho_k)^{3} e^{-(\log \rho_k)^{3/2}}.
 \end{equation}
Thanks to \eqref{eq:V12prima}, \eqref{eq:V22prima},
\eqref{eq:stimaderivataprima} and \eqref{eq:equivEk3d} 
we have that for some $C>0$ and for big $k$
\begin{equation}\label{eq:collectme1}
\begin{split}
  \abs{\mathbb V(x)} \leq  & C \abs{ \chi'_k(r) n_k^{\frac74}} 
   \leq  C \rho_{k}^{-\frac{1-\delta}{2}} 
   (\log\rho_k)^{3} e^{-(\log \rho_k)^{3/2}} 
   n_k^\frac74  \\
   \leq & C (\log\rho_k)^{3} 
   \rho_{k}^{-\frac{1-\delta}{2}}   
   \leq C (\log \abs{x})^3 \abs{x}^{-\epsilon},
\end{split}
\end{equation}
using \eqref{eq:relazioneepsdelta} in the last inequality.

If $s \geq (\log \rho_k)^{-3/2}$, we have
\begin{equation}\label{eq:V12seconda}
  \abs{\V_{21}(x)} \leq  \frac{\abs{\chi_k'(r)} \abs{E_{k+1}(x)} \abs{\chi_k(r)} \abs{E_{k+1}(x)}}%
  {\abs{\chi_k(r) E_{k+1}}^2 + \abs{E_{k}(x)}^2} 
  \leq  \frac{\abs{\chi_k'(r)} }{\abs{\chi_k(r)}}
\end{equation}
and
\begin{equation}\label{eq:V22seconda}
  \abs{\V_{22}(x)} \leq  \frac{\abs{\chi_k'(r)} \abs{E_{k+1}(x)} \abs{E_{k}(x)}}%
  {\abs{\chi_k(r) E_{k+1}}^2 + \abs{E_{k}(x)}^2} 
  \leq  \frac{\abs{\chi_k'(r)} }{2\abs{\chi_k(r)}}.
\end{equation}
If $s \in [(\log \rho_k)^{-3/2}, \frac12]$,  thanks to
\eqref{eq:V12seconda}, \eqref{eq:V22seconda}, \eqref{eq:stimaderivata} and the explicit expression for $\chi_k$, for 
large $k$ we have that
\begin{equation}\label{eq:collectme2}
  \abs{\mathbb V (x)} = O(\rho_{k}^{-\frac{1-\delta}{2}}) s^{-2} 
  \leq C (\log \rho_k)^3 \rho_{k}^{-\frac{1-\delta}{2}} 
  \leq  C (\log\abs{x})^3 \abs{x}^{-\epsilon},
\end{equation}
while in the case $s \in (1/2,1)$, since $\chi$ is non-decreasing, we
have that for some $C>0$
\begin{equation}\label{eq:collectme3}
  \abs{\mathbb V (x)} \leq C \rho_{k}^{-\frac{1-\delta}{2}}
  \leq  C \abs{x}^{-\epsilon}.
\end{equation}
Gathering \eqref{eq:collectme1}, \eqref{eq:collectme2} and \eqref{eq:collectme3}
we have \eqref{eq:decayV2d} for all $r \in [\rho_{k1},\rho_{k2}]$.
 Finally, the case $r \in [\rho_{k2},\rho_{k3}]$ is analogous and will
 be omitted in this proof.

\subsection{Conclusion of the proof}
The remaining part of the proof is analogous to the proof in the case
that $n=2$:
the construction of $u$ in $\{ \abs{x} \leq \rho_{k_0} \}$ is done as in \Cref{sec:defn.origin},
the construction of $u$ in $\{ \abs{x} \geq \rho_{k_0} \}$ as in \Cref{sec:defn.glueing} 
and the study of the decay of $u$ at infinity as in \Cref{sec:decay.u}.

\printbibliography

\end{document}